\documentclass[12pt]{amsart}

\usepackage{amscd}
\usepackage{mathabx}
\input bookman.sty
\boldmath

\usepackage{comment}
\excludecomment{excl}

\usepackage{xcolor}

\usepackage{helvet}

\usepackage{graphics}
\usepackage{amssymb}
\usepackage{amsxtra}
\usepackage{amsmath}
\usepackage{mathrsfs}

\usepackage[arrow, matrix]{xy}
\xyoption{frame}

\theoremstyle{plain}

\newtheorem{theo}{Theorem}[section]
\newtheorem{lemm}[theo]{Lemma}
\newtheorem{prop}[theo]{Proposition}

\newtheorem{coro}[theo]{Corollary}

\theoremstyle{definition}

\newtheorem{defi}[theo]{Definition}
\newtheorem{rema}[theo]{Remark}

\newfont{\rmm}{cmr10 scaled 1000}

\newfont{\itt}{cmsl10 scaled 1000}

\newfont{\rM}{cmr10 scaled 1700}

\newcounter{lemma}[section]

\newcounter{tempcounter}

\newcommand{\lb}{\label}

\newcommand{\rrf}[1]{(\ref{#1})}

\renewcommand{\a}{\alpha}
\renewcommand{\b}{\beta}
\newcommand{\g}{\gamma}

\renewcommand{\t}{\theta}
\renewcommand{\l}{\lambda}

\newcommand{\m}{\mu}

\renewcommand{\r}{\rho}

\newcommand{\s}{\sigma}

\renewcommand{\o}{\omega}

\renewcommand{\Re}{\text{\rm Re }}
\newcommand{\D}{\Delta}

\renewcommand{\O}{\Omega}

\renewcommand{\AA}{{\mathcal A}}
\newcommand{\BB}{{\mathcal B}}

\newcommand{\DD}{{\mathcal D}}
\newcommand{\EE}{{\mathcal E}}
\newcommand{\FF}{{\mathcal F}}

\newcommand{\II}{{\mathcal I}}
\newcommand{\JJ}{{\mathcal J}}

\newcommand{\LL}{{\mathcal L}}

\newcommand{\OO}{{\mathcal O}}

\newcommand{\cc}{{\mathbb{C}}}

\newcommand{\hh}{{\mathbb{H}}}

\newcommand{\qq}{{\mathbb{Q}}}
\newcommand{\rr}{{\mathbb{R}}}

\newcommand{\ttt}{{\mathbb{T}}}

\newcommand{\zz}{{\mathbb{Z}}}

\newcommand{\LLLL}{{\mathscr{L}}}

\newcommand{\OOOO}{{\mathscr{O}}}
\newcommand{\PPPP}{{\mathscr{P}}}

\newcommand{\id}{\text{id}}

\renewcommand{\Im}{\text{\rm Im }}

\newcommand{\Id}{\text{\rm Id}}

\newcommand{\GL}{\text{\rm GL}}

\newcommand{\bere}{\begin{rema}}
\newcommand{\bede}{\begin{defi}}

\renewcommand{\beth}{\begin{theo}}
\newcommand{\bele}{\begin{lemm}}
\newcommand{\bepr}{\begin{prop}}
\newcommand{\beeq}{\begin{equation}}
\newcommand{\bega}{\begin{gather}}
\newcommand{\begaa}{\begin{gather*}}
\newcommand{\been}{\begin{enumerate}}

\newcommand{\bedee}{\begin{defii}}
\newcommand{\bethh}{\begin{theoo}}
\newcommand{\belee}{\begin{lemmm}}
\newcommand{\beprr}{\begin{propp}}

\newcommand{\beco}{\begin{coro}}

\newcommand{\beal}{\begin{aligned}}

\newcommand{\enre}{\end{rema}}

\newcommand{\enco}{\end{coro}}
\newcommand{\enpr}{\end{prop}}
\newcommand{\enth}{\end{theo}}
\newcommand{\enle}{\end{lemm}}
\newcommand{\enen}{\end{enumerate}}
\newcommand{\enga}{\end{gather}}
\newcommand{\engaa}{\end{gather*}}
\newcommand{\eneq}{\end{equation}}
\newcommand{\enal}{\end{aligned}}

\newcommand{\bq}{\begin{equation}}
\newcommand{\bqq}{\begin{equation*}}

\renewcommand{\leq}{\leqslant}
\renewcommand{\geq}{\geqslant}

\newcommand{\ove}{\overline}

\newcommand{\sbs}{\subset}

\newcommand{\sut}{~such~that~}

\newcommand{\wrt}{with respect to}
\newcommand{\ho}{homomorphism}

\newcommand{\ma}{manifold}

\newcommand{\Prf}{{\it Proof.\quad}}

\newcommand{\chart}{\Phi_p:U_p\to B^n(0,r_p)}
\newcommand{\atlas}{\{\Phi_p:U_p\to B^n(0,r_p)\}_{p\in S(f)}}

\newcommand{\pr}{\partial}

\newcommand{\qt}{\hfill\triangle}
\newcommand{\qs}{\hfill\square}

\newcommand{\pa}{\vskip0.1in}

\newcommand{\arrh}[3]
{
\xymatrix{
{#1} \ar[r]^<<<<{#2}  &{#3}
}
}

\newcommand{\arrr}[1]
{\arrh {}{#1}{}}

\newcommand{\arrto}
{\xymatrix{{} \ar@{|-{>}}[r]  & {} } }

\newcommand{\arrinto}
{\xymatrix{{} \ar@{^{(}->}[r]  & {} } }

\newcommand{\OT}{Oeljeklaus\ and\ Toma}
\newcommand{\ot}{Oeljeklaus-Toma}

\begin{document}

\title
[On generalized OT manifolds]
{On a generalization of Inoue and  Oeljeklaus-Toma manifolds}
\author{Hisaaki Endo  and  Andrei Pajitnov}
\address{Department of Mathematics
Tokyo Institute of Technology
2-12-1 Ookayama, Meguro-ku
Tokyo 152-8551
Japan}
\email{endo@math.titech.ac.jp}
\address{Laboratoire Math\'ematiques Jean Leray 
UMR 6629,
Universit\'e de Nantes,
Facult\'e des Sciences,
2, rue de la Houssini\`ere,
44072, Nantes, Cedex}                    
\email{andrei.pajitnov@univ-nantes.fr}

\thanks{} 
\begin{abstract}
In this paper we construct a family of complex analytic
manifolds that  generalize Inoue surfaces and 
Oeljeklaus-Toma manifolds. 
To a  matrix $M$ in $SL(N,\mathbb{Z})$
%having $s$ real  simple eigenvalues 
satisfying some mild conditions
on its characteristic polynomial 
we associate a manifold $T(M,\DD)$ 
(depending on  an auxiliary parameter $\DD$).
This manifold fibers over $\ttt^s$ 
with fiber $\ttt^{N}$ (here $s$ is the number of 
real eigenvalues of $M$);
the monodromy matrices are certain 
polynomials of  the matrix
$M$.
The basic difference of our construction 
from the preceding ones is that we admit non-diagonalizable 
matrices $M$ and the monodromy of the above fibration 
can also be non-diagonalizable. We prove that 
for a large class of non-diagonalizable matrices
$M$ the manifold $T(M,\DD)$ 
does not admit any K\"ahler structure and 
is not 
homeomorphic to any of Oeljeklaus-Toma manifolds.
\end{abstract}
\keywords{}
\subjclass[2010]{32J18, 32J27, 57R99}
\maketitle
\tableofcontents
\newcommand{\ssll}{\mathrm{SL}}

\section{Introduction}
\label{s:intro}
\subsection{Background}
\label{su:bckgr}

In 1972 M. Inoue \cite{Inoue1974}
constructed complex surfaces
having remarkable properties:
they have second Betti number equal to zero
and contain no complex curves.
Inoue surfaces
attracted a lot of attention.
It was proved by F. Bogomolov \cite{Bogomolov}
(see also the works of J. Li, S.-T. Yau, and F. Zheng 
\cite{LYZh} and  \cite{LYZ},
and A. Teleman \cite{T}) 
that each complex surface of class $VII_0$ 
with $b_2(X)=0$ and containing no complex curves 
is isomorphic to an Inoue surface.
Inoue surfaces are not algebraic, and moreover 
they do not admit K\"ahler metric.

Let us say that a matrix $M\in \ssll(2k+1,\zz)$ 
is of  type $\II$, if it has only one real 
eigenvalue which is irrational and strictly greater than $1$.
Inoue's construction associates to every 
matrix $M\in \ssll(3,\zz)$ of  type  $\II$
a complex surface $S_M$ obtained as a quotient 
of $\hh\times \cc$ by action of a discrete group
(here $\hh$ is the upper half-plane).
This manifold fibers over  $S^1$ with fiber $\ttt^3$
%and the monodromy of this fibration equals 
%the diffeomorphism of $\ttt^{3}$ determined by 
%$M^{\top}$.

Inoue's construction was generalized
to higher dimensions in the work of 
%G. K. Sankaran
%\cite{Sankaran}
 K. Oeljeklaus and M. Toma 
\cite{OT}.
The construction of Oeljeklaus and Toma
uses algebraic number theory.
It starts with an algebraic field $K$ of degree $n$.
Such field admits $n$ embeddings to $\cc$;
denote by $s$ the number of real embeddings;
then $n=s+2t$.
\OT~ construct an action of a certain semidirect 
product $\zz^s \ltimes \zz^{2t+s}$
on $\hh^s\times \cc^t$,
the quotient is a compact complex manifold
of complex dimension $s+t$.
It has interesting geometric properties,
in particular, it does not admit a K\"ahler metric.
The original Inoue surface corresponds to the 
algebraic number field generated by the eigenvalues
of the matrix $M$.

The Oeljeklaus-Toma  manifolds 
({\it OT-manifolds } for short)
have very interesting geometric properties, studied in \cite{OT};
in particular, they do not admit K\"ahler metric.
These manifolds were recently studied by many authors.
In the work of L. Ornea, M. Verbitsky, and V. Vuletescu 
\cite{OVV}
it is shown that in many cases the OT-manifolds do not contain 
proper analytic subvarieties. 
In the article 
\cite{IO} of
N. Istrati and A. Otiman
the De Rham cohomology of OT-manifolds is computed. 
The paper of D. Angella, M. Parton, and V. Vuletescu
\cite{APV}
is dedicated to the proof of the rigidity of the complex structure 
of the OT-manifolds. 
The non-existence of complex curves in OT-manifolds is proved in the paper 
\cite{SV} of S. Verbitsky.

In the paper \cite{EP}
we introduced another generalization of the Inoue surfaces.
Namely for every matrix $M$ of type $\II$ 
we constructed a complex non-K\"ahler \ma~ $T_M$ 
that fibers over a circle and the monodromy of the fibration 
equals $M^\top$.
If the matrix $M$ is non-diagonalizable,
then the \ma~ $T_M$ is not homeomorphic to any of the manifolds 
constructed by Oeljeklaus and Toma.

\subsection{Outline of the paper}
\label{su:outline}

In the present paper we generalize our construction  
to  the case of matrices $M$ with several real 
eigenvalues. Let $M\in SL(N,\zz)$. 
%be a matrix
%having  at least one real and at least one non-real eigenvalue.
Assuming  some  restrictions on the matrix $M$ 
(see Definition \ref{d:propJ})
we construct
a  \ma~ $T(M,\DD)$
that fibers over
$\ttt^s$ (where $s$ is the number of real eigenvalues of $M$)
with fiber $\ttt^{s+2n}$ and the monodromy group 
isomorphic to $\zz^s$.
The monodromy transformations are linear automorphisms of the fiber
deduced from the matrix $M$ (see Theorem \ref{t:inoue_s}).

The contents of the paper section by section is as follows.
In Section \ref{s:twi-act} we discuss some 
special properly discontinuous
actions of semidirect product $\zz^s\ltimes \zz^k$ 
on 
$
(\rr^*_+)^s
\times
\rr^k$.
We call them {\it twisted-diagonal actions}.
In Section \ref{s:twi-matrices} we show how to construct
a twisted-diagonal actions
from a matrix
$M$ of special type.
Section \ref{s:dir-fam}
 is the central technical part of the work,
here we study the {\it Dirichlet families} of
integer polynomials associated to a matrix.
The construction of the  generalized \ot~  manifold 
is the subject of  the main theorem of the 
paper (Theorem \ref{t:inoue_s})
which is stated and proved in Section \ref{s:inoue_s}.
This theorem  associates to each matrix $M$ 
of type $\JJ$ (see Definition \ref{d:propJ})
and a Dirichlet family $\DD$ for $M$ 
a complex manifold $T(M,\DD)$.
Sections \ref{s:multi-tor}  -- \ref{s:Tm-OT}
are about properties of the manifolds $T(M,\DD)$.
We concentrate on the case when the matrix $M$ is 
non-diagonalizable, since this case is the most different
from the manifolds construsted in \cite{Inoue1974}, \cite{OT}.
We prove in particular that if $M$ is non-diagonalizable,
and the family $\DD$ is {\it primary} (see Definition \ref{d:primaryDD})
then:
\pa
-- The manifold $T(M,\DD)$ does not admit any LCK-structure,
in particular it does not admit a structure of a K\"ahler 
manifold. 
\pa
-- The manifold $T(M,\DD)$ is not homeomorphic to any 
of OT-manifolds.
\pa 
In section \ref{s:Tm-OT} we  show that some of 
OT-manifolds are biholomorphic to manifolds $T(M,\DD)$.

%\newpage

\section{Twisted diagonal actions}
\label{s:twi-act}

In this subsection we consider some special groups of diffeomorphisms
of the space
$
(\rr^*_+)^s
\times 
E
$
where $E$ is a finite-dimensional real vector space.

Let $k=\dim E$.
Let $\BB$ be any basis in $E$ 
and denote by $\LL$ the free abelian
group generated by $\BB$, so that $\LL\approx \zz^k$.
The group $\LL$ acts on $E$ by translations.
Let 
$\FF=(\xi_1,\ldots,   \xi_r )$
be any family of elements of
$(\rr^*_+)^s$. Each of vectors $\xi_i$ 
can be considered as a diagonal matrix in $GL(s,\rr)$.
These matrices commute with each other and so we obtain 
an action of $\zz^r$ on $\rr^s$. We can restrict it to 
$(\rr^*_+)^s$.
The direct product of the two actions above 
gives an action of the abelian group $G_0=\zz^r\times \LL$ on 
$(\rr^*_+)^s\times E$. This action  will be referred to as 
{\it diagonal action} of $G_0$.

Now we shall generalize this notion so as to define 
an action of a certain semidirect product of $\zz^r$ and $\LL$
on $(\rr^*_+)^s\times E$.

Let
$\AA=\{A_1, \ldots, A_r\}\sbs \GL(E)$
be a family of commuting linear automorphisms.
Assume that every $A_i$ preserve the lattice $\LL$,
and moreover, determines an automorphism of $\LL$. 
Define an action of $\zz^r$ on $\LL$ by the following rule:
the $i$-th free generator $\tau_i$ of $\zz^r$
acts as $A_i$. 
Consider the corresponding semidirect product
$
G=
\zz^r
\ltimes
\LL.
$

\bede\lb{d:twi-di-action}
Define
an action of $G$ 
on 
$(\rr^*_+)^s\times E$ as follows.

1) The subgroup $\LL$ of $G$ acts by translations on $E$.

2) The generator $\tau_i$ of the group $\zz^r$ 
acts as follows
$$
\tau_i(x_1, \ldots, x_s, v) =
(\xi_i^{(1)}x_1, \ldots, \xi_i^{(s)}x_s, A_i(v))
$$
This action will be called 
{\it twisted diagonal action }
associated to $\FF, \LL$ and $\AA$.
The  quotient of $(\rr^*_+)^s\times E$ by this action will be 
 denoted by 
$Y(\FF,\LL,\AA)$.
The action 
 will be called 
{\it special}
if at least one of the linear
maps $A_i:E\to E$
has no eigenvalue 1.
\end{defi}

\bede\lb{d:log-basis}
The map
$$
(\rr^*_+)^s
\to
\rr^s,
\ \ \ 
(x_1, \ldots, x_s )
\arrto
(\log x_1, \ldots, \log  x_s )
$$
will be denoted by $\lg$.
For a family $\FF=(\xi_1,\ldots,   \xi_r )$
of elements of $
(\rr^*_+)^s$
we denote by $\lg\FF$
the family 
$(\lg \xi_1, \ldots, \lg  \xi_r )$
of elements of $\rr^s$.
A family $\FF$
is called {\it a log-basis}
if $s=r$ and $\lg\FF$ is a basis in $\rr^s$.
$\qt$
\end{defi}

\bede\lb{d:non-deg-twi-action}
The twisted diagonal action 
associated to $\FF, \LL$ and $\AA$
will be called 
{\it non-degenerate}
if $\FF$ is a log-basis of
$(\rr^*_+)^s$.
\end{defi}

\bepr\label{p:twidiag-compact}
Any non-degenerate twisted diagonal action
is properly discontinuous and the quotient
is a compact \ma.
\enpr
\Prf
Consider the homeomorphism
$$
\psi
:
(\rr^*_+)^s
\times  E
\arrr {\lg\times \id} 
\rr^s\times E
$$
The action of $G$ 
on 
$
(\rr^*_+)^s
\times  E
$
induces an action of $G$ on 
$
\rr^s\times E
$
which satisfies
the following:

1) The lattice $\LL$ acts
by translations on $E$.

2) The action of the generator 
$\tau_i\in\zz^s$ on $\rr^s\times E$ 
is given by the formula
$$
\tau_i(a,b)
=
(a+\lg \xi_i, A_i(b)).
$$
The vectors 
$
(\lg \xi_1, \ldots, \lg  \xi_s )
$
form a basis in $\rr^s$ and 
the proposition follows. $\qs$

Denote the matrix of $A_i$ in the basis $\BB$ by $\mu_i$,
then $\mu_i\in GL(k,\zz)$.

\bepr\lb{p:quotient}

\been
\item 
There is a fibration 
$p:Y(\FF,\LL, \AA)\to \ttt^s$
with fiber $\ttt^k$
and monodromy matrices equal to $\mu_1, \ldots, \mu_s$.
\item 
This fibration admits a cross-section $\s$.
\enen 
\enpr
\Prf 
To prove 1) define a diagonal  action of $\zz^s$ on $(\rr^*_+)^s$
by the formula 
$$
\tau_i(x_1, \ldots, x_s) =
(\xi_i^{(1)}x_1, \ldots, \xi_i^{(s)}x_s)
$$
and observe that 
the first projection 
$(\rr^*_+)^s\times E
\to
(\rr^*_+)^s$
commutes with the actions
of $\zz^s$ and determines a map
$Y(\FF,\LL,\AA)\to \ttt^s$
which is clearly a fibration 
with fiber $E/\LL\approx \ttt^k$
and  monodromy maps $\mu_i$ as required. 
As for the point 2) 
the cross-section
is defined by the 
inclusion 
$$
(\rr^*_+)^s
\approx
(\rr^*_+)^s\times \{0\} 
\sbs 
(\rr^*_+)^s \times E.
$$
$\qs$

%\bepr\lb{p:compl-structure}
%Assume that $E$ is a complex vector space and the 
%linear maps $A_i$ preserve the complex structure.
%Then the manifold 
%$Y(\AA, \FF,\LL)$
%has a natural structure of a complex manifold 
%of dimension $s+k$.
%\enpr
%\Prf
%The 

%The point 1) follows from the fact that $G$ acts by biholomorphisms.

%\newpage

\section{Twisted diagonal actions associated with an integer 
matrix}
\label{s:twi-matrices}

\bede\lb{d:propertyJ0}
We say that a polynomial $P\in\zz[t]$ is of {\it type $\JJ_0$}
if it has at least one real root and 
at least one imaginary root and 
all real roots are simple.

Let $M\in SL(N,\zz)$.  We say that $M$ is a matrix of {\it  type $\JJ_0$}
if its characteristic polynomial if of type $\JJ_0$.
 $\qt$
\end{defi}
\noindent
In this section we outline a procedure of constructing from 
a matrix 
$M$ of type $\JJ_0$ a twisted diagonal
action. This procedure requires some further assumptions 
on $M$ which will be discussed in subsequent sections.
Denote by $\a_1, \ldots, \a_s$ 
the real eigenvalues of $M$ and by
$\b_1, \ldots, \b_k, \bar \b_1,  \ldots, \bar \b_k$
the imaginary ones. We have then $0<s, 0<k$. The size of the matrix 
equals $s+2n$ with $n\geq k$.
%We assume that $\a_s>1$ for all $s$.
%In this section we associate to $M$ a twisted diagonal
%action.

 %\newpage
 Denote by $f_M:\cc^{s+2n}\to \cc^{s+2n}$
 the linear map corresponding to $M$.
We have a decomposition of $\cc^{s+2n}$
into a direct sum of complex $f_M$-invariant subspaces
$$
\cc^{s+2n} = V\oplus W\oplus \bar W
$$
where $V$ is generated by the eigenvectors 
$
a_1, \ldots , a_s\in \rr^{s+2n}
$ 
of $M$, corresponding to the real eigenvalues 
$
\a_1, \ldots, \a_s
$, and $W$ is the generalized eigenspace corresponding to 
all the imaginary eigenvalues.
Let 
$b_1, \ldots , b_n$
be a basis of $W$.
Consider the vectors 
\bq\lb{f:def-v}
v_i=\Big( a^{(i)}_1, \ldots , a^{(i)}_s, b^{(i)}_1,
\ldots, b^{(i)}_n \Big) \in \rr^s\times \cc^n
\approx \rr^{s+2n}.
\end{equation}
(here $1\leq i \leq s+2n$)
%and $w_i\in\hh$).

\bele\lb{l:basis}
The family $\BB=(v_1, \ldots , v_{s+2n})$ 
is a basis of $\rr^{s+2n}$.
\enle
\Prf
It suffices to prove that 
the determinant of the matrix $Q$ formed 
by the columns of the coordinates
of $(s+2n)$ real vectors 
$a_1, \ldots, a_s, \Re b_1,\ \Im b_1, 
\ldots ,\Re b_n,\ \Im b_n$
is non-zero.
This determinant equals to a non-zero multiple 
of the determinant of the following
family of complex vectors:
$(a_1, \ldots, a_s, b_1, \bar b_1, \ldots b_n, \bar b_n)$,
and this last family is a basis. $\qs$ 

The next definition is basic for our purposes.
\bede\lb{d:R-and-U}
Denote by  $R$  the matrix of the restriction of the linear 
map $f_M$ to $W$. 
Denote by $\D$ the diagonal $(s\times s)$-matrix
with diagonal entries 
$\a_1, \ldots, \a_s$. Let $U\in GL(s+2n,\cc)$
be the direct sum $\D\oplus R$ of these matrices.
Define  a linear automorphism $g_0$
of $\cc^s\times \cc^n$  by the following formula
\bq\lb{f:def-g0}
g_0:\cc^s\times \cc^n\to \cc^s\times \cc^n; \ \ \  (w_1, \ldots, w_s, z) \mapsto 
(\a_1w_1, \ldots, \a_sw_s, R^\top z)
\end{equation}
(here $w_i\in\cc$ for every $i$ and $z\in \cc^n$).
Its matrix equals  $U^\top$.
Observe that the above formula determines also an 
automorphism of $\rr^s\times \cc^n$, and a biholomorphism of
$\hh^s\times \cc^n$
onto itself. These two isomorphisms will be denoted by the same 
symbol $g_0$ by a certain abuse of notation.
 $\qt$
\end{defi}

We will need a lemma from linear algebra.

\bele\lb{l:change_bas}
Let $A$ be an $(m\times m)$-matrix. We will 
denote the corresponding mapping
$\cc^m\to\cc^m$ by the same letter.
Let $L$ be an $A$-invariant subspace,
and $\BB=(b_1, \ldots b_k)$ a basis in $L$.
Denote by $B$ the $(m\times k)$-matrix formed by 
the columns of the vectors of $\BB$, and by $R$ 
the matrix of $A~|~L:L\to L$ in the basis $\BB$. 
Then 
\bq\lb{f:ch-b}
AB=BR.
\end{equation}
\enle
\Prf
Denote by $\EE$ the canonical basis of $\cc^m$. Observe
that $AB$ is the matrix of the composition 
$L\arrinto \cc^m\arrr A \cc^m$ \wrt~ the bases 
$\BB$ and $\EE$. 
Further, $BR$ 
is the matrix of the composition 
$L\arrr A L\arrinto \cc^m$ \wrt~ the  same bases.
The lemma follows. $\qs$

\bere\lb{r:ch-bss}
If $L=\cc^m$, the formula 
\rrf{f:ch-b}
is equivalent
to the classical formula
$R=B^{-1}AB$, 
relating the matrices
of a linear map in two different bases. 
\enre

We will now apply the preceding Lemma 
to the basis  $\BB$ from Lemma \ref{l:basis}.

\bepr\lb{p:lattice-invar}
The lattice $\LL$ generated by $\BB$ 
is invariant \wrt~ $g_0$,
and the matrix of the map $g_0$ in $\BB$
equals $M^\top$.
\enpr
\Prf 
Put $u_i=( a^{(i)}_1, \ldots , a^{(i)}_s)\in\rr^s$
and
$
w_i
= (b^{(i)}_1,
\ldots, b^{(i)}_n )\in \cc^n$,
so that $v_i=(u_i,w_i)$;
here $1\leq i \leq s+2n$.
%Denote by $\D$ the diagonal $(s\times s)$-matrix
%with diagonal entries 
%$\a_i$. 
We have
$(\D u_i)_j=\a_j a^{(i)}_j$;
this equals 
$\sum_k M_{ik} a^{(k)}_j$
since
$a_i$ is the eigenvector of $M$
with eigenvalue $\a_i$. Therefore 
$\D u_i=\sum_k M_{ik} u_k$. 

Let $B$ be the matrix formed by the columns of coordinates of vectors 
$b_i$. Then the vector $R^\top w_i$ is the $i$-th column 
of the matrix $R^\top B^\top$. This 
matrix equals $(MB)^\top$
by Lemma \ref{l:change_bas}
and $(MB)_{ij} = \sum_k M_{ik} b_k^{(j)}$.
Therefore $g_0(v_i)=\sum_k M_{ik} v_k
=\sum_k (M^\top)_{ki} v_k
$ and 
the proposition is proved. $\qs$
\pa
This proposition allows us to construct 
a twisted diagonal action of $\zz\ltimes \LL
\approx \zz\ltimes\zz^{s+2n}$
on $(\rr^*_+)^s\times \Big(V\times W\Big)$,
following the procedure from Section \ref{s:twi-act}.
The input data is $(\FF, \LL, \AA)$
where $\FF$ consists of just one vector 
$(\a_1,\ldots, \a_s)$
and 
$\AA$ consists of just one matrix $U^\top$.
For the case $s=1$ this action is nothing else than 
the one constructed in our previous paper \cite{EP}.

For the case $s>1$ this action 
is not  non-degenerate.
To construct a non-degenerate action we 
need a family of $s$ matrices, commuting with each other.
We will choose these matrices to be of the form $D_i(U^\top)$
where $D_i$ are integer polynomials of special type,
which will be described now.

\newcommand{\veca}{\overrightarrow {\a}}
\newcommand{\vecg}{\overrightarrow {\g}}

%Let $D_1, \ldots , D_s\in \zz[t]$
For any 
polynomial $D\in \zz[t]$
the vectors 
$a_1, \ldots , a_s\in\rr^{s+2n}$
are eigenvectors of $D(M)$ 
with  eigenvalues 
$D(\a_1), 
\ldots ,
D(\a_s).
$
Put
$${\veca}= (\a_1, \ldots, \a_s), \ \ 
D({\veca})
=
\Big(D(\a_1), \ldots, D(\a_s)\Big).
$$

\bede\lb{d:diri-fam}
A family 
$
\DD
=(D_1, \ldots ,D_s)
$
of integer polynomials
will be called {\it a Dirichlet family for $M$} or 
{\it a $D$-family } for short,
if the following three conditions hold:

$\D 1)$\ $D_i(M)\in SL(s+2n,\zz)$ for every $i$;

$\D 2)$\ $D_i(\a_j) > 0$ for every $i, j$;

$\D 3)$\ The family 
$
\DD(\veca) = \Big(D_1(\veca), \ldots ,D_s(\veca)\Big)
$
is a log-basis of 
$(\rr^*_+)^s.
$
\end{defi}

\bere\lb{r:positivity}
Observe that if a family $\DD=(D_1, \ldots ,D_s)
$ 
satisfies only the conditions $\D 1)$ and $\D 3)$
it is easy to construct from it a Dirichlet family,
namely, put $\tilde \DD=(\tilde D_1, \ldots ,\tilde D_s)$
where
$\tilde D_i=D^2_i$.
\enre

We will construct such families later, now we are going to explain how
to construct a twisted-diagonal action from a $D$-family.
For every $i$ with $1\leq i \leq s$
the automorphism 
$$g_i=D_i(g_0) :\rr^s\times \cc^n
\to
\rr^s\times \cc^n;
$$
preserves the lattice $\LL$;
its  matrix in the canonical basis  equals $D_i(U^\top)$.
and in basis $\BB$ it equals $D_i(M^\top)$.

Applying Definition \ref{d:twi-di-action}
to the input data $(\FF, \LL, \DD)$
where $\FF=(D_1(\veca), \ldots ,D_s(\veca))$
and $\DD=(D_1(U^\top),\ldots, D_s(U^\top))$
we 
obtain a 
twisted-diagonal 
action 
of
$G=\zz^s\ltimes \zz^{s+2n}$
on 
 $
(\rr^*_+)^s 
 \times 
 E$ where 
$E= \rr^s\times \cc^n.
 $ 
We have a natural diffeomorphism
$$
\psi : (\rr^*_+)^s \times \Big(\rr^s\times \cc^n\Big)
\approx \hh^s\times \cc^n.
$$
\noindent
In the coordinates $\psi$ the action of $g_i$
%$i$-th generator of $\zz^s$
is written as follows:
\bq\lb{f:def-act}
g_i(w,z) 
=
 \Big(D_i(\a_1)w_1,\ldots, D_i(\a_s)w_s,
 D_i(R^\top)z
 \Big)
 \end{equation}
 (where $w\in \hh^s, \ z\in \cc^n$).
 Thus every $g_i$ is a biholomorphism 
 of $\hh^s\times \cc^n$ and we obtain the 
 following.
 
 \bepr\lb{p:complex}
 Let $\DD$ be a Dirichlet family of polynomials for a matrix $M$
 of type $\JJ_0$.
 Then the formula 
 \rrf{f:def-act}
 determines a twisted-diagonal action
 of the corresponding semidirect product
 $G=\zz^s\ltimes \zz^{s+2n}$
 on $\hh^s\times \cc^n$, and we have the following:
 
\been\item
The quotient  $Y$ has a natural structure of a complex
\ma~ of dimension $s+n$.
\item There is a fibration $p:Y\to \ttt^s$ with fiber $\ttt^{s+2n}$;
the action of generators of the monodromy group $\zz^s$
is given by the matrices 
$D_1(M^\top), \ldots, D_s(M^\top)$.
\item The fibration $p$ admits a cross-section.
\enen 
\enpr
\Prf
The part 1) follows from the fact that $G$
acts by biholomorphisms,
together with Proposition  \ref{p:twidiag-compact}.
The parts 2) and 3) 
follow from  Proposition 
\ref{p:quotient}. $\qs$

 \bere\lb{r:s-equals1}
For $s=1$ the family consisting of one polynomial $D_1(t)=t$
is obviously a Dirichlet family for any matrix $M$ of type $\JJ_0$.
Thus we recover here  the construction of generalized Inoue manifolds
from  \cite{EP}, \S 2.
\end{rema} 
 \pa

 %\newpage 
 
 \section{Constructing Dirichlet families}
 \lb{s:dir-fam}
 
 In this section we construct Dirichlet families 
 associated to a given matrix $M\in SL(N,\zz)$.
 The main result of this section is  Theorem \ref{t:constr-diri}.
 The principal algebraic tool in this construction is Proposition \ref{p:diri},
 based on the Dirichlet's unit theorem.

 \subsection{A corollary of the Dirchlet's unit theorem}
 \lb{su:coro_diri}

\bepr
\lb{p:diri}
%Let 
%$M\in SL(N,\zz)$
%denote by $Q$ its characteristic polynomial.%
Let $P\in\zz[t]$ be 
an irreducible polynomial.
Assume that $P$ has at least one real root, and at least one 
imaginary root. Denote the real roots of $P$ by
$\g_1, \ldots, \g_l$, and let $2q$ be the number of imaginary roots.
%so that $l+2q=N$.
%Let $e_1, \ldots, e_l$ be the corresponding real eigenvectors.
Put $\vecg=(\g_1, \ldots, \g_l)$.

Then there are polynomials $P_1, \ldots, P_l\in\zz[t]$ \sut:

1) $P_i$ are invertible in $\zz[t]/P$.

2) For every $i, j$ we have 
$P_i(\g_j)>0$.

3) The vectors $P_i(\vecg)=(P_i(\g_1), \ldots, P_i(\g_l))$
with  $1\leq i \leq l$
 form a log-basis in $(\rr^*_+)^l$.
 
4) The properties 1) -- 3) hold if we replace $P_i$ by 
any other polynomial $\tilde P_i$ \sut~ $\tilde P_i \equiv P_i \mod P$.
\enpr
 \Prf 
 Consider the algebraic extension $K$ of $\qq$ 
 corresponding to $P$, that is, $K=\qq[t]/P$.  Let $m=\deg P$, then $m=l+2q=(K:\qq)$.
  An embedding $K \hookrightarrow \cc$ 
is called {\it real} if its image is     in $\rr$;
an embedding which is not real is called {\it complex}.
%(\cite{BS}).
The field $K$ admits $l$  real embeddings
and  $2q$  complex embeddings.
Denote the real embeddings by 
$\s_1, \ldots, \s_l$ and the complex  embeddings
by $\s_{l+1}, \ldots, \s_{l+2q}$;
we can assume             that 
$\s_i=\ove{\s_{q+i}}$ for $i\geq l+1$.
The map 
$$\s:K\to\rr^l\times \cc^q;
\ \ \ 
\s(x)=\big(\s_1(x), \ldots , \s_{l+q}(x)\big)
$$
is an embedding (known as {\it geometric representation}
of the field $K$, see \cite{BS}, Ch. II, \S 3).
The map $\LLLL:K^*\to\rr^{l+q}$ defined as follows:
$$
\LLLL(x)=
$$
$$\Big(\log|\s_1(x)|, \ldots , 
\log|\s_l(x)|,
2\log|\s_{l+1}(x)|,
\ldots,
2\log|\s_{l+q}(x)|
\Big)
$$
is called {\it the logarithmic representation map}.
 
 The ring $ \zz[t]/P$ is an order in $K$, denote it by $\OO$.
  Dirichlet's unit theorem 
 (see \cite{BS}, Ch. II, \S 4.3) guarantees 
 the existence of $l+q-1$ elements $u_1, \ldots, u_{l+q-1}$
 in $\OO$, which are invertible in $\OO$ 
 and generate a free abelian group in $\OO^*$.
 Moreover, these units can be chosen so as their 
 $\LLLL$-images form a basis in the hyperplane 
 $$
 H=\Big\{\l\in \rr^{l+q} ~|~
 \l_1+...+\l_{l+q}=0\Big\},
 $$
 see \cite{BS}, Chapter II, \S 4.3.
 Denote by $p$ the projection of $\rr^{l+q}$ 
 onto the space $\rr^l$ corresponding to the first $l$ coordinates.
  Since $q\geq 1 $ we can assume (re-numbering the 
units $u_i$ if necessary), that the vectors
  $p(\LLLL(u_1)), ... ,p(\LLLL(u_l))$
form a basis in $\rr^l$.
 
 Pick any polynomials $P_i\in\zz[t]$ (here $1\leq i \leq l$)
 \sut~ $P_i(t)=u_i$ in $\OO$.
 Every real embedding of $K$ sends $t$ to one of the real roots of 
 $P$, so 
renumbering the roots $\g_i$ if necessary we can assume that 
$\s_i(t)=\g_i$.
Therefore 
$$\s_j(u_i)=\s_j(P_i(t))=P_i(\g_j).$$
Replacing the units $u_i$ by their squares 
 $u_i^2$ if necessary 
 we can assume that 
 $P_i(\g_j)>0$.
The property  2)  of the statement of the 
Proposition follow.
As for 3) it suffices to observe that the vector $\lg P_i(\vecg)$
equals 
$p(\LLLL(u_i))$.
The  point 4) is obvious.
 $\qs$

 \beco\lb{c:irreducible-case}
 Let $M\in SL(N,\zz)$;
 assume that the  characteristic polynomial $P$ of $M$
 satisfies the assumptions of Proposition 
 \ref{p:diri}. Then there exists a Dirichlet family
 of polynomials for $M$.
 \enco
 \Prf 
 %The real roots of $P$ are simple since $P$ is irreducible,
 %so $M$ is of type $\JJ_0$. 
 The family 
 $(P_1,\ldots, P_l)$
 of polynomials   from Proposition \ref{p:diri}
 satisfy  the properties 
 $\D 2)$ and $\D 3) $ (they follow from 
 the properties 2) and 3) of Proposition \ref{p:diri}).
  Further, observe that every matrix $P_i(M)$ is invertible.
 Indeed there exists a polynomial $T_i$ \sut~ 
 $P_iT_i\equiv 1\ \mod (P)$. Therefore $P_i(M)T_i(M)=Id$
 and $P_i(M)$ is invertible. 
 
 It remains only to show that $\det P_i(M) >0$.
 This property follows from the next  Lemma,
 the proof of which is an easy consequence of the Jordan decomposition theorem 
 and will be omitted. 

 \bele\lb{l:det-posi}
 Let $M$ be a square matrix with entries in $\rr$, and $P\in\rr[t]$.
 Assume that for every real eigenvalue $\a$ of $M$ 
 we have $P(\a)\geq 0$. Then $\det P(M) \geq 0$. $\qs$
 \enle

 $\qs$

 \subsection{Dirichlet families: a particular case}
 \lb{su:constru-diri1}
 
In this subsection we  begin the construction of  Dirichlet families;
we consider here 
a particular class of matrices.
% a matrix $M\in SL(N,\zz)$.
We will use the following terminology for polynomials
with integer coefficients.

\bede\lb{d:coprimes}
Let $A,B\in\zz[t]$. We say that $A$ and $B$ are {\it coprime}
if their greatest common divisor is $1$.
We say that $A$ and $B$ are {\it strongly coprime}
if the sum of the principal ideals $(A)$ and $(B)$
equals $\zz[t]$.
\end{defi}

\bede\lb{d:type-J1}
An integer polynomial $C$ is said to be of type $\JJ_1$
if $C=B_0\cdot B_1$ with $B_0, B_1\in\zz[t]$, 
and
\been\item 
$B_0$ has no real roots.
\item 
$B_1$ is irreducible, it 
 has at least one real root and 
 at least one imaginary root. 
\item
%The ideals generated by 
$B_0$ and $B_1$ are strongly coprime.
%, that is,
% $(B_0)+(B_1)=\zz[t]$.
\enen 

We say that a matrix $M\in SL(N,\zz)$
is of type $\JJ_1$ if its characteristic polynomial
is of type $\JJ_1$. $\qt$
 
 \end{defi}

 % ex:twidiag-tor

 \bere\lb{r:jord}
 This condition allows  $B_0$ 
 to have multiple complex roots.
 \enre

  \bepr
  \lb{p:diri-fam2}
  If $M$ is of type $\JJ_1$, 
  then there exists a Dirichlet family 
  of polynomials for $M$.
  \enpr
  \Prf
  Denote by $C$ the characteristic polynomial of $M$.
  The ring $\zz[t]ex:twidiag-tor
/C$ has natural projections
  $
  \phi_0:
  \zz[t]/C
  \to 
  \zz[t]/B_0$
  and
  $\phi_1:
  \zz[t]/C
  \to 
  \zz[t]/B_1$.
  The product of these projections 
  $$\phi=
  \phi_0
  \times 
  \phi_1
  :
  \zz[t]/C
  \to 
  \zz[t]/B_0
  \times 
  \zz[t]/B_1$$
  is an isomorphism
  by the 
  Chinese remainder theorem (see \cite{H}, Th 2.25, p. 131).
  Denote the real roots of $B_1$ by $\a_1, \ldots, \a_s$.
  Apply Proposition \ref{p:diri}
  to $B_1$;
  we obtain polynomials $P_1, \ldots , P_s\in\zz[t]$. For $1\leq i \leq s$
  let $Q_i\in\zz[t]$
  be any polynomial invertible $\mod B_0$, and 
  $D_i\in \zz[t]$
  be a polynomial \sut~
  \begin{equation}\lb{f:*}
     D_i\equiv Q_i \mod ( B_0),{\rm \ \ \ \ and \ \ \ \ }
  D_i\equiv P_i \mod( B_1 ).
  \end{equation}
  We claim that the family 
  $\DD=(D_1, \ldots , D_s)$
  is a $\DD$-family for $M$.
  
  Since $D_i(\a_j)=P_i(\a_j)$, the formula \rrf{f:*}
  allows us to deduce the property $\D 2)$ from property 
  2) of Proposition \ref{p:diri}.
  As for the property $\D 3)$ it follows immediately from property 3) 
  of \ref{p:diri}.
  
  The polynomial $D_i$ 
  is invertible in $\zz[t]/C$ since $\phi$ is an isomorphism
  and both
  $\phi_0(D_i)$
  and
  $\phi_1(D_i)$
  are invertible.
  Therefore $D_i(M)$ is invertible over $\zz$.
  Applying Lemma 
  \ref{l:det-posi}
  we deduce that $\det D_i(M) =1$.

  %An argument similar to the proof of 
%Corollary \ref{c:irreducible-case}
 % shows that $\det D_i(M)=1$.
  
   $\qs$
  
  \bere\lb{r:primaryD}
  There is a large choice for polynomials $Q_i$.
  In particular the polynomials $Q_i(t)=t$ will do.
    Any Dirichlet family corresponding to this case 
  will be called {\it primary}.
  \end{rema}
  
  Here is one example of a polynomial of type $\JJ_1$.

 \bepr\lb{p:exampleJ1}
 Let $B_1(t)=t^4+3t^3+3t^2+3t+1, A(t)=t^2+1$.
 Let $q$ be any natural number $\geq 1$, and put $B_0(t)=(A(t))^q$.
 Then the polynomial $C(t)=B_0(t)B_1(t)$ is of type $\JJ_1$.
 \enpr
 \Prf
 It is easy to check that $B_1$ is irreducible over $\zz$
 and has two real roots (both negative),
 and two imaginary roots.
 The resultant $R(B_1, A)$ equals $1$,
 therefore 
  the ideal generated by $B_0$ and $A$ contains $1$,
  so that $B_0$ and $A$ are strongly coprime.
 We have also
 $R(B_1, A^q)=(R(B_1, A))^q=1$, therefore
 $B_0(t)$ and $B_1(t)$ are coprime. 
 $\qs$ 
%  equals $\zz[t]$.

 Using recent results  from the theory of irreducible
 polynomials one can generalize the previous
  example to construct polynomials of type $\JJ_1$ 
  of higher degrees. Here is one more example.
  \bepr\lb{p:exampleJ1bis}
 Let $B_1(t)=t^{12}+8t^8+20t^7+19t^6+20t^5+8t^4+1$, and $A(t)=t^2+1$ as before.
  Let $q$ be any natural number $\geq 1$, and put $B_0(t)=(A(t))^q$.
 Then the polynomial $C(t)=B_0(t)B_1(t)$ is of type $\JJ_1$.
 \enpr
  \Prf 
  The coefficients of $B_1$ are all positive,
  the degree of $B_1$ is less than 31,
  and the natural number $B_1(10)$ is  prime. Therefore 
  by a theorem of M. Filaseta and  S. Gross
  (\cite{FG}, Th. 1) the polynomial $B_1$ is irreducible.
  In particular, its roots are simple.
  We have $ B_1(-1)=-3$, therefore 
  $B_1$ has at least  two real roots. The polynomial $f'(x)$ 
  is divisible by $x^3$, and can have at most 9 real roots, 
  therefore $B_1$ has at least  one imaginary root.
    The resultant $R(B_1, B_0)$ equals 1, similarly to the 
  Proposition \ref{p:exampleJ1}, and the proof is over. $\qs$

 \subsection{Dirichlet families: general  case}
 \lb{su:constru-diri2}
 
 In this subsection we construct Dirichlet families
 in more general situation than 
 in the previous subsection.
 
 \bede\lb{d:propJ}
 We say that 
an integer polynomial $C$ is {\it of type $\JJ$}
if 
$C=B_0B_1\cdots B_k$
 where $B_i$
 have the following properties:
 
 $\JJ$1) $B_0$ has no real roots.
 
 $\JJ$2) For every $i\geq 1$ 
 the polynomial $B_i$ is irreducible,  
 it has at least one real root and 
 at least one imaginary root.
  
 $\JJ$3) 
 For every $i,j \geq 0$ 
 the  polynomials $B_i, B_j$
 are  strongly coprime.
 \pa
 
 We say that a matrix $M\in SL(N,\zz)$
 is {\it of type $\JJ$}
 if its characteristic polynomial
 is of type  $\JJ$. $\qt$ 
 \end{defi} 
 
 We have obviously $\JJ_0\Leftarrow \JJ \Leftarrow \JJ_1 $.
   Here is one example of a polynomial of type $\JJ$ which is not of type $\JJ_1$.

 \bepr\lb{p:exampleJ}
 Let   $A(t)=t^2+1$ and 
 $$B_1(t)=t^4+3t^3+3t^2+3t+1, \ \ 
 B_2(t)=t^4+4t^3+3t^2+4t+1.$$
 Let $q$ be any natural number $\geq 1$, and put $B_0(t)=(A(t))^q$.
 Then the polynomial $C(t)=B_0(t)B_1(t)B_2(t)$ is of type $\JJ$.
 \enpr
 \Prf
 It is easy to check that  each of the polynomials
 $B_1, B_2$ is irreducible over $\zz$
 and has two real roots (both negative),
 and two imaginary roots. As for the resultants of these polynomials
 we have
 $$
 R(B_1, A)=1,
 \ \ 
  R(B_2, B_1)=1,
  \ \ 
  R(B_2, A)=1;
  $$
  the proposition follows. 
  $\qs$

 \beth\lb{t:constr-diri}
 Let $M\in SL(N,\zz)$
 be a matrix 
 of type $\JJ$.
  Then there exists a $D$-family  of 
 polynomials for $M$.
 \enth 
 \Prf
 Denote by $s_j$ the number of real roots of $B_j$,
 let $s=\sum_{j=1}^k s_j$.
  The ring $\zz[t]/C$ has natural projections 
 $\phi_j:\zz[t]/C \to \zz[t]/B_j$.
 Consider the product of these projections
 \bq\lb{f:chin-rem}
 \phi: \zz[t]/C \to \zz[t]/B_0 \times 
 \zz[t]/B_1\times ...  \times \zz[t]/B_k.
 \end{equation}
 Since the polynomials $B_j$ are 
 strongly pairwise coprime, $\phi$ is an isomorphisms
 by the Chinese remainder theorem (see \cite{H}, 
 Th 2.25, p. 131).

For $1\leq j \leq k$ apply Proposition \ref{p:diri} to the 
polynomial $B_j$; we obtain polynomials $P_{j,i}$,
where $1\leq i \leq s_j$. 
The polynomials of the Dirichlet family in construction
will be indexed by couples $(j,i)$ with  $1\leq j\leq k$
and $1\leq i \leq s_j$.
Let $(j,i)$  be such couple. 

\pa
Let $E_{j,i}\in\zz[t]$
be any polynomial 
invertible $\mod B_0$.
%(here $1\leq j \leq k$).
%For $j\geq 1$
\pa
Let $D_{j,i}$ be an integer polynomial \sut~

%\phi_0\Big(\ove%For an integer polynomial $A$ we will denote by $\ove A$ 
%its image in $\zz[t]/C$.
$${D_{j,i}}\equiv E_{j,i}\ {\rm mod} (B_0), $$
$$
D_{j,i} \equiv P_{j,i}\ {\rm mod} (B_j), $$
$$D_{j,i}
\equiv 1\ {\rm mod} (B_\mu) 
\ \ 
{\rm if \ \ }
\mu\not= j\ \ 
{\rm and  \ \ }\mu\not= 0.
$$
We can choose such polynomials since the map 
 \rrf{f:chin-rem} is an isomorphism.
%We obtain a family $\{D_{j,i}\}$ 
%( where 
%$1\leq j \leq k$ 
%and $1\leq i \leq s_j$)
%containing $s$ polynomials.
We claim that the family $\DD=\{D_{j,i}\}$ 
is a $D$-family for $M$. 

Let us first prove the  
property $\D 2)$.
Denote by 
$\a_{j,i}$
the real roots of $B_j$
(here $1\leq i \leq s_j$).
%(here $j \geq 1$ and $1\leq i \leq s_j$).
We have $P_{j,i}(\a_{j,\l})>0$ for every $\l$;
this implies $
D_{j,i}(\a_{j,\l}) >0 
$
since
$D_{j,i}
\equiv
P_{j,i}\ {\rm mod} (B_j)
$
and 
$
B_j(\a_{j,\l})=0$.
Similarly, 
if $\mu\not= j$, then 
$D_{j,i}(\a_{\mu,\l})=1$
for every $\l$ and 
the property $\D 2)$ holds.

Let us proceed to the property $\D 1)$.
For an integer polynomial $A$ we will denote by $\ove A$ 
its image in $\zz[t]/C$.
Observe that
for every $\m\geq 0 $
the element
$\phi_\mu(\ove{D_{j,i}})
$
is invertible, therefore 
$\ove{D_{j,i}}
$
is also invertible, since $\phi$ is an isomorphism.
This implies that the matrix $D_{j,i}(M)$ is invertible
over $\zz$.
Applying Lemma 
  \ref{l:det-posi}
  we deduce that $\det D_{j,i}(M) =1$.

To prove $\D 3)$  observe that the $s\times s$ matrix 
$\log D_{j,i}(\a_{\mu,\l})$
is block-diagonal.
It has $k$ blocks of sizes 
$s_1 , \ldots ,s_k$
and the $j$-th block 
equals  the     $(s_j\times s_j)$-matrix
$\log D_{j,*}(\a_{j,*})$
with $1\leq * \leq s_j$.
The determinants of these matrices
are non zero by the property 3) of Proposition
\ref{p:diri}. The proof of Theorem \ref{t:constr-diri}
is now over. $\qs$

\begin{defi}\lb{d:primaryDD}
There is a large choice for polynomials $E_{j,i}$.
  In particular the polynomials $E_{j,i}(t)=t$ will do.
    Any Dirichlet family corresponding to this case 
  will be called {\it primary}.
  \end{defi}

\section{The manifold $T(M,\DD)$}
\lb{s:inoue_s}

Now we summarize
the contents of the previous section
to obtain the generalized OT-manifold. In this section 
$M$ denotes a matrix in $SL(N,\zz)$ of type $\JJ$.
This matrix has 
$s$ simple  real eigenvalues with $s>0$  and 
 at least one imaginary eigenvalue, and we have  $N=s+2n$ with $n>0$.
 Let $\DD=(D_1, \ldots , D_s)$
 be any D-family for $M$.
 Such families exist by Theorem \ref{t:constr-diri}
 of Subsection \ref{su:constru-diri2}.
 The family $\DD(M)=(D_1(M), \ldots , D_s(M))$
  of integer matrices determines an action of
  $\zz^s$ on $\zz^{s+2n}$.
  Denote by $G_\DD$ the group
  associated to  this action, so that 
  $G_\DD\approx \zz^s\ltimes \zz^{s+2n}$.
  The next result 
is the main theorem of the 
present paper. It follows immediately from 
Proposition
\ref{p:complex}.
  
  \beth\lb{t:inoue_s}
  The group 
$G_\DD$
acts by biholomorphisms on the space
 $\hh^s\times \cc^n$.
 This action is properly discontinuous and the quotient $T(M, \DD)$
 is a compact complex $(n+s)$-dimensional manifold.
 We have a natural fibration $T(M, \DD)\to \ttt^s$
 with fiber $\ttt^{s+2n}$, and the monodromy group isomorphic
 to $\zz^s$. The action of generators of the monodromy group
 on $\ttt^{s+2n}$ is given by 
 the matrices $D_1(M^\top), \ldots , D_s(M^\top)$.
 $\qs$
 \enth 

 %\bede\lb{d:gen-in}
 %The manifold constructed in the previous theorem will be denoted
 %by $T(M, \DD)$.
 %\end{defi}
 
 %\bere\lb{r:s-equals1}
%For $s=1$ the family consisting of one polynomial $D_1(t)=t$
%is obviously a Dirichlet family for any matrix $M$ of type $\JJ_0$.
%Thus we recover here  the construction of generalized Inoue manifolds
%from  \cite{EP}, \S 2.
%\end{rema} 
% \pa

 \subsection{Relation with generalized Inoue manifolds from \cite{EP}}
 \lb{su:rel-gen-inoue}
 
 In this subsection we explain the relation of manifolds $T(M,\DD)$
 with the generalized Inoue \ma s, constructed in \cite{EP}.
 Recall from \cite{EP}
 that a matrix $M\in SL(2n+1,\zz)$ is called {\it matrix of type $\II$}
 if it has only one real eigenvalue which is simple, positive
 and irrational.
 In \cite{EP} we associated to each such  
 matrix a complex non-K\"ahler 
 manifold $T_M$ of complex dimension $n+1$
 endowed with a fibration $T_M\to S^1$, whose fiber is a torus
 $\ttt^{2n+1}$.

 A matrix of type $\II$ is clearly a matrix of type $\JJ_0$
 with $s=1$. In this case the family consisting of 
 one polynomial $D_1(t)=t$
is  a Dirichlet family.
  The twisted-diagonal action of the semidirect product
$\zz\ltimes \zz^{2n+1}$ 
corresponding to $D_1$
is the same as the action constructed in \cite{EP}
and the manifold 
$T(M,\DD)$ 
is biholomorphic to the manifold $T_M$ from \cite{EP}.

%\newpage

\section{Mapping multi-tori
and their homological properties}
\lb{s:multi-tor}

From the topological  point of view
Inoue surfaces are 
just the mapping tori of the corresponding 
linear map of $\ttt^{3}$.
In this section 
we discuss a generalization of the mapping torus
construction, which leads to a simple 
description  of the topological type of $T(M,\DD)$
and enables us to establish its basic topological
properties. 

\bede\lb{d:multi-t}
Let $X$ be a topological space and
$\phi:\zz^s\to Homeo(X) $ a group \ho.
Consider the diagonal action of $\zz^s$
on $\rr^s\times X$
($\zz^s$ acting on $\rr^s$ by translations).
The {\it mapping multi-torus of $\phi$}
is the quotient of 
$\rr^s\times X$
by this action.
It is denoted by 
$MT(\phi)$.
\end{defi}

\bere\lb{r:multitor}
\been\item We have a natural fibration 
$p_\phi~:~MT(\phi)\to\ttt^s$ with fiber $X$.
\item 
If $s=1$ then the space 
$MT(\phi)$
is the mapping torus of the homeomorphism
$\phi(1):X\to X$.
\item
Choose a basis $t_1,\ldots, t_s$ in $\zz^s$,
let $G_l$ be the free abelian subgroup generated
by 
$t_1, \ldots, t_l$.
Let $\phi_k=\phi~|~G_k$.
The map $\phi(t_{k+1})$ determines a 
self-homeomorphism
$\psi_{k+1}: MT(\phi_k) \to MT(\phi_k)$
and 
$MT(\phi_{k+1})$
is the mapping torus
of $\psi_{k+1}$.
\enen
\enre

\bede\lb{d:special}
A homomorphism $\ \phi:\zz^s\to Homeo(X) $ 
is called {\it special}
if 
\been\item 
The space $X$ is connected.
\item 
All the maps $\phi(g)$ for $g\in \zz^s$
have a common fixed point.
\item 
There is $g\in \zz^s$
\sut~ 
the homomorphism
$(\phi(g))_*-\Id:
H_1(X,\qq)
\to
H_1(X,\qq)
$
is injective.
\enen
\end{defi}

Observe that the property 2) from the above definition
implies existence of a cross-section
of the fibration $p_\phi$.

\bepr\lb{p:betti_spec}
If  $\phi$
is special, then
 $b_1(MT(\phi))=s$,
and the fibration $p_\phi$
induces an isomorphism in $H_1(\ - \ , \qq)$.
\enpr
\Prf
Choose a basis $t_1,\ldots, t_s\in\zz^s$
in such a way
that $1$ is not an eigenvalue
of the \ho~ $(\phi(t_1))_*$
of $H_1(X,\qq)$.
Let $B_l=MT(\phi_l)$.
We will now prove two following assertions
by induction in $l$ with $l\geq 1$.

\begin{quote}
1) \ \ $
b_1(B_l, \qq)=l$, 
\ \ and  

\noindent
2) $(\phi(t_{l+1}))_* : 
H_1(B_l, \qq)
\to
H_1(B_l, \qq)$
equals $\Id$.
\end{quote}

Assuming that these assertions hold
for a given value of $l\geq 1 $
write the Milnor's exact sequence (see \cite{Milnor})
$$
H_1(B_l)
\xrightarrow {(\phi(t_{l+1}))_* -\Id}
H_1(B_l)
\to
H_1(B_{l+1})
\to 
H_0(B_l)
\to
H_0(B_l).
$$

The space $B_l$ is connected, 
therefore $b_1(B_{l+1})\leq l+1$.
The fibration $B_{l+1}\to \ttt^{l+1}$
has a cross-section, therefore 
$b_1(B_{l+1})\geq l+1$.
The part 1) of the induction step
follows. The part 2)  follows from the fact that 
the cross-section induces isomorphism in $H_1$.
The induction step is now complete. 
A similar argument proves the assertion for 
$l=1$, here we use the fact that 
1 is not an eigenvalue of the homomorphism
$(\phi_1)_*$.
$\qs$

\bepr\lb{p:fund-gr}
Assume that the properties
1) and 2)
of Definition 
\ref{d:special}
hold for a group \ho~
$\phi:\zz^s\to Homeo(X) $.
Then the group 
$\pi_1(MT(\phi))$
is the semidirect product
$\zz^s\ltimes \pi_1(X)$
\wrt~ the corresponding
\ho~ \ 
$\tilde\phi:\zz^s\to Aut(\pi_1(X))$.
(Namely,  $\tilde\phi$ 
sends any element $g\in\zz^s$ to
the automorphism in $\pi_1$ induced by $\phi(g)$.)
\enpr
\Prf
Follows immediately 
from the exact sequence
of the fibration $p_\phi$. $\qs$

\bede\lb{d:prod-diag}
 Let 
 $\r:\zz^s\to GL(k,\zz)$
 be a \ho.
 Denote by $G_{s,k}(\r)$
 the corresponding
 semidirect 
 product
 $\zz^s\ltimes \zz^k$.
 We say that $\r$
 is {\it special}
 if there is  
 $g\in \zz^s$ \sut~ $1$ is not an eigenvalue
 of $\r(g)$.
 We say that $\r$
 is {\it diagonalizable}
 if for every $g\in \zz^s$ the 
 matrix $\r(g)$
 is diagonalizable over $\cc$.
 \end{defi}
 
 Each \ho~ $\r:\zz^s\to GL(k,\zz)$
 determines a \ho~ 
 $\bar\rho:\zz^s\to Homeo(\ttt^k)$;
 if $\r$ is special then $\bar\r$ is also special.
 
 \beco\lb{c:bettiG}
 If $\r:\zz^s\to GL(k,\zz)$
 is special
 then 
 $b_1( G_{s,k}(\r))=s$.
 $\qs$
 \enco

 \bele\lb{l:two}
 Assume that a \ho~ $\r:\zz^s\to GL(k,\zz)$
 is special.  Let 
 $f_1, f_2:G_{s,k}(\r)\to \zz^s$
 be epimorphisms.
 Then there is an isomorphism
 $\psi:\zz^s\to \zz^s$ 
 \sut~ 
 $f_2= f_1\circ\psi$.
 \enle 
 \Prf
 Consider the Hurewicz homomorphism modulo torsion: 
 \bq\lb{f:hur}
 h: G_{s,k}(\r)\to H_1(G_{s,k}(\r))/Tors\approx \zz^s.
 \end{equation}
 The \ho~ $f_1$ factors through 
 $h$, so that
 $f_1=\psi_1\circ h$.
 The map $\psi_1$ is then an epimorphism
 of $\zz^s$ onto itself, therefore it is an isomorphism.
 Similarly $f_2=\psi_2\circ h$
 and the proof is over. $\qs$
 
 \bepr\lb{p:isos}
  Assume that  $G_{s,k}(\r) \approx G_{t,m}(\t)$
 where $\r$ is
 special.
 Assume that $b_1( G_{t,m}(\t)) = s$. 
  Then $s=t,\  k=m$ and there is an isomorphism 
 $\psi:\zz^s\to\zz^s$ \sut~ 
 the representations 
 $\r$ and $\t\circ \psi$ are equivalent. 
 In particular, if $\r$ is diagonalizable, then 
 $\t$ is also diagonalizable.
 \enpr
 \Prf 
 %We have $b_1(G_{s,k}(\r))=s=b_1(G_{t,m}(\t))=t$.
 Let $\Psi:G_{s,k}(\r)\to G_{t,m}(\t)$ be an isomorphism.
 Denote by $j: G_{t,m}\to \zz^s$ 
 the Hurewicz homomorphism modulo torsion
 (similarly to \rrf{f:hur}).
 Lemma \ref{l:two} 
 implies that there is an isomorphism
 $\psi$ \sut~ the following diagram is commutative.
  $$\xymatrix{ 
  0 \ar[r] &\zz^k \ar[d]^{\psi_0}\ar[r] &G_{s,k}(\r) 
  \ar[d]^{\Psi}\ar[r]^h &\zz^s  \ar[d]^{\psi}\ar[r] &0 \\
  0 \ar[r] &\zz^m \ar[r] &G_{t,m}(\t) \ar[r]^j &\zz^s  \ar[r] &0
 }$$
 Then $\psi_0$ is necessarily an isomorphism,
 therefore $m=k$ and for each $g\in\zz^s$ the matrix
 $\r(g)$ is conjugate to $\t(\psi(g))$ via the isomorphism 
 $\psi_0$.
 $\qs$
 
 Let us now return to generalized Inoue manifolds.
 Let $M\in SL(s+2n,\zz)$ and $\DD$ be a D-family for $M$.
 Every map
 $g_i=D_i(g_0): 
 \rr^s\times \cc^n
 \to
  \rr^s\times \cc^n
  $
  preserves the lattice
  $\LL$
  and determines a self-homeomorphism
  $\bar g_i$ of the torus
  $X
  =
   \rr^s\times \cc^n
   /
   \LL.
   $
Define a map 
$
\phi:\zz^s\to 
Homeo(X)
$
by $\phi(t_i)=g_i$;
     then 
   the manifold $T(M,\DD)$
   is clearly homeomorphic to the 
   mapping multi-torus
   $
   MT(\phi).
   $
   
   \beco\lb{c:pi1-dir}
   The group $\pi_1(T(M,\DD))$
   is isomorphic to the semidirect product
   $\zz^s\ltimes \zz^{s+2n}$
   where the action of $\zz^s$ on $\zz^{s+2n}\approx \LL$
   is given by the homomorphism
   $t_i\mapsto D_i(M^\top)$.
   \enco 
   Denote by $G(\DD)$
   the  subgroup of $SL(s+2n,\zz)$
   generated by $D_i(M^\top)$ (where $1\leq i \leq s$).

 \bede\lb{d:specialDir}
 The Dirichlet family
 $\DD$ is called
 {\it special}, if in 
 the subgroup $G(\DD)$ there is at least
 one matrix
 without
 eigenvalue $1$.
 \end{defi}
 For example
 in the case $s=1$ the Dirichlet
 family consisting of one polynomial $D_1(t)=t$
 is special.
 \bepr\lb{p:specD}
 Let $M$ be a matrix of type $\JJ$
 with characteristic polynomial $C=B_0\cdot B_1\cdot \ldots \cdot B_k$
 (see Definition \ref{d:propJ}).
 Then any primary Dirichlet family
 is special.
 (In particular this is the case when $B_0=1$.)
 \enpr
 \Prf
 The matrix
 $$
 \Big( \log D_i(\a_j) \Big)_{1\leq i,j \leq s}
 $$
 has a non-zero determinant.
 Therefore some linear combination
 of its lines has only strictly positive coordinates.
 Thus for some $n_1, \ldots , n_s\in\zz$ 
 we have
 $$
 \sum_{i=1}^s n_i\log D_i(\a_j)\not=0 
{\rm \
 and \
 }
 \prod_{i=1}^s D_i^{n_i}(\a_j) \not=1
 \ \ 
 {\rm for \ every }
 \ \ 
 j.
 $$
 Let $N=\prod_{i=1}^s D^{n_i}_i(M^\top)$.
 Then $N\in G(\DD)$ and 
 $1$ is not an eigenvalue of
 $N~|~ V$.
 Further, $
 N~|~ (W\oplus \bar W)
 =
 M^\top~|~ (W\oplus \bar W)
 $
 since the family $\DD$ is primary.
 Therefore
 all eigenvalues of $N$ 
 in $(W\oplus \bar W)$
 are imaginary. $\qs$
 
 \beco\lb{c:b-1}
 Let $\DD=(D_1, \ldots , D_s)$ 
 be a special Dirichlet
 family for $M$. Then $b_1(T(M,\DD))=s$.
 $\qs$
 \enco

%\newpage

 \section{LCK structures on $T(M,\DD)$}
 \lb{s:lck}
 
 In this section we show that if $M$ is 
  non-diagonalizable, and 
 $\DD$ is a primary 
 Dirchlet family 
 for $M$, then the manifold $T(M,\DD)$
 does not admit a structure of a locally conformally K\"ahler manifold
 (in particular it is not K\"ahler). The argument follows basically the 
 lines of the proof of the corresponding property of \OT~ manifolds 
 (\cite{OT}, Proposition 2.9), although we do not need the condition 
 $s=1$ from  \cite{OT}. 
 
 \beth\lb{t:non-lck}
 Let $M$ be a matrix of type $\JJ$, and $\DD$ be a primary 
 Dirchlet family for $M$. Assume that at least one of the 
 matrices $D_i(M)$ is 
 non-diagonalizable. Then the manifold $T(M,\DD)$
 does not admit a structure of a locally conformally K\"ahler manifold.
 \end{theo}
 \Prf
 Denot $T(M,\DD)$ by $Y$ for brevity. 
 We have a tower of coverings
 $$
Y'= \hh^s\times\cc^n\arrr \chi \bar Y \arrr \pi Y
$$
where $\pi$ is induced from the universal covering $\rr^s\to\ttt^s$
via the fibration  $Y\to \ttt^s$ from Theorem \ref{t:inoue_s},
and the covering $\chi$ is the universal covering
for $\bar Y$.
Assume that $Y$ admits an LCK-metric, let $\O$ be the corresponding $2$-form
on $Y$
so that $d\O=\o\wedge\O$, and $d\o=0$.
Then $d(\pi^*\O)=\pi^*\o\wedge\pi^*\O$.
The map $\pi$ induces the sero \ho~ in $H_1(\cdot, \rr)$,
therefore we have $\pi^*\o=df$.
The form $\t=e^{-f}\pi^*\O$ on
$\bar Y$ has the following properties:
\pa
A) $\t$ is closed.

B) For every $g\in \zz^s$
we have $g^*\t=\r(g)\cdot \t$, where $\r(g)=e^{-\langle [\o], g\rangle}$.

C) The bilinear form  $\tilde \t$ associated to $\t$ 
via the following formula:
\bq\lb{f:scal-pr}
\tilde \t(x)(h,k)=\t(x)(h,\sqrt {-1} k)
\end{equation}
is positive definite. 
\pa

The manifold $\bar Y$ is diffeomorphic to $\ttt^{s+2n}\times (\rr^*_+)^s$,
and so there is a natural action of the group $\ttt^{s+2n}$ on it.
Denote by $\t_1$ the 2-form 
obtained by
averaging the form $\t$ \wrt~ this action 
and the usual Haar measure on 
$\ttt^{s+2n}$.
It is clear that $\t_1$ 
satisfies 
the  conditions A) and C) above.
The next lemma shows that the condition B) is also true for $\t_1$.

\bele\lb{l:condB}
For every $g\in \zz^s$
we have $g^*(\t_1)=\r(g)\cdot \t_1$.
\enle
\Prf
Let $\xi\in 
\ttt^{s+2n}$.
Observe that 
$$
g^*(\xi^*\t)
=
g^*\xi^*(g^*)^{-1}g^*\t
=
\r(g) 
g^*\xi^*(g^*)^{-1}
\t.$$
Therefore
$$
\int_{\ttt^{s+2n}} 
g^*(\xi^*\t) d\xi
=
\r(g)
\int_{\ttt^{s+2n}}
g^*\xi^*(g^*)^{-1}\t d\xi.
$$
The diffeomorphism
$\xi\mapsto g\xi\g^{-1}$
of 
$
\ttt^{s+2n}$
onto itself preserves  volumes, 
since all matrices $D_i(M)$ are in $SL(2n+s,\zz)$.
The   changement of variables theorem
implies that 
$$\int_{\ttt^{s+2n}} (g^{-1}\xi g)^*\t d\xi 
=
\int_{\ttt^{s+2n}} \t d\xi 
$$
and we have 
$g^*\t_1 = \r(g)\t_1$. $\qs$

The  restriction of $\t_1$ to any subspace
$\ttt^{s+2n}\times v$
where $v\in  (\rr^*_+)^s$
has constant coefficients \wrt~ the 
usual basis of
$\O^2(\ttt^{s+2n})$. 
Moreover let $\l=\chi^*(\t_1)$,
then the coefficients of $\l$
\wrt~ the usual basis of
$\O^2(\hh^s\times\cc^n)$
do not depend on $\vec z$
(where $(\vec w, \vec z)$ denote the coordinates on 
$\hh^s\times\cc^n$).
Thus the restriction of $\l$ to any subspace $w\times \cc^n$
writes as
$$
\sum_{i,j} \l_{ij}(\vec w)dz_i\wedge  d\bar z_j.
$$
\bele\lb{l:lambda}
The coefficients $\l_{ij}$
do not depend on $\vec w$.
\enle
\Prf 
Write 
$$
\l
=
\sum_{i,j} g_{ij}(\vec w)dw_i\wedge d\bar w_j
+
\sum_{i,j} h_{ij}(\vec w)dw_i\wedge d\bar z_j
$$
$$
+
\sum_{i,j} k_{ij}(\vec w)dz_i\wedge d\bar w_j
+
\sum_{i,j} \l_{ij}(\vec w)dz_i\wedge d\bar z_j.
$$
The coefficients of $d\l$
at
$dw_k\wedge dz_i\wedge d\bar z_j$
and
at
$d\bar w_k\wedge dz_i\wedge d\bar z_j$

equal
$ \frac {\pr \l_{ij}(w)}{\pr w_k}$, and 
respectively
$ \frac {\pr \l_{ij}(w)}{\pr \bar w_k}$.
Since $d\l=0$, the lemma follows. $\qs$

Put $\l_0=\l~|~(1\times \cc^n)$.
The symmetric bilinear form $\tilde \l_0$
associated to $\l_0$ via the formula 
\rrf{f:scal-pr}
is a scalar product on $\cc^n\approx \rr^{2n}$.
Since the family $\DD$ is primary, the
restriction of any matrix $D_i(M)$ to $\cc^n$ 
equals $R^\top$ and is 
therefore non-diagonalizable.
Put $U=D_i(M)$.
Then we have 
$\mu(U\xi, U\eta)=C\mu(\xi,\eta)$
for $\xi, \eta \in \cc^n$, so that 
the matrix $U_0=U/\sqrt {C}$ preserves the scalar product.
This contradicts to the assumption that $U$ is non-diagonalizable. $\qs$

 %\newpage

\section{Relations with the Oeljeklaus-Toma construction}
\lb{s:Tm-OT}

%\subsection{Comparison of manifolds $T_M$ with OT-manifolds 
%} \lb{su:compar}

This section is about relations of the manifolds $T(M,\DD)$
with the manifolds constructed by  \OT~ in \cite{OT}
({\it OT-manifolds} for short).
We begin by a brief recollection of 
Oeljeklaus-Toma construction, 
then we show that some of 
OT-manifolds are biholomorphic to manifolds $T(M,\DD)$.
Then we show that 
in the non-diagonalizable case 
the manifold 
$T(M,\DD)$
is not homeomorphic to any of OT-manifolds.

\subsection{Construction of  OT-manifolds }
\lb{su:constrOT}

%Let us first recall the construction from \cite{OT}
%(in a slightly modified  terminology).
Let $K$ be an algebraic number field.
An embedding $K \hookrightarrow \cc$ 
is called {\it real} if its image is     in $\rr$;
an embedding which is not real is called {\it complex}.
%(\cite{BS}).
Denote by $s$ the number  of real embeddings 
and by $n$ the number  of complex embeddings.
Then $(K:\qq)=s+2n$.
Let $\s_1, \ldots, \s_s$ be the real embeddings
and $\s_{s+1}, \ldots, \s_{s+2n}$
be the complex  embeddings;
we can assume             that 
$\s_i=\ove{\s_{n+i}}$ for $i\geq s+1$.
%Put $r=s+t$.
The map 
$$\s:K\to\rr^s\times \cc^n;
\ \ \ 
\s(x)=\big(\s_1(x), \ldots , \s_{s+n}(x)\big)
$$
is an embedding (known as {\it geometric representation}
of the field $K$, see \cite{BS}, Ch. II, \S 3).
Let $\OOOO$ be any order in $K$, then 
$\s(\OOOO)$ is a full lattice in 
$
\rr^s\times \cc^n.
$
Denote by $\OOOO^*$ the group of all units of $\OOOO$.
The Dirichlet Unit Theorem (see \cite{BS}, Ch. II, \S 4 , Th. 5)
says that the group 
$\OOOO^*/Tors$
is a free abelian group of rank $s+n-1$.
Assume that $n\geq 1$.
Choose any elements  $u_1, \ldots, u_s$ 
of $\OOOO^*$ generating 
in $\OOOO^*/Tors$
a free abelian subgroup
of rank $s$. A unit $\l\in\OOOO$ 
will be called {\it positive}
if $\s_i(\l)>0$ for every $i\leq s$.
Replacing
$u_i$ by $u_i^2$ if necessary 
we can assume that 
every $u_i$ is positive.
The subgroup $U$ of 
$\OOOO^*$
generated by $u_1, \ldots, u_s$ 
acts on $\OOOO$
and we can form 
the semidirect product $\PPPP=U\ltimes \OOOO$.
The group $\PPPP$
acts on 
$\cc^r=\cc^s\times \cc^n$
as follows:
\begin{itemize}
 \item 
 any element $\xi\in\OOOO$
 acts by translation by vector
 $\s(\xi)\in \rr^s\times \cc^r.$
 \item any element $\l\in U$ acts
 as follows:
 $$
\l\cdot (z_1, \ldots, z_{s+n})
=
(\s_1(\l)z_1, \ldots, \s_{s+n} (\l)z_{s+n}).
$$
\end{itemize}
For $i\leq s$ the numbers $\s_i(\l) $
are  real and positive, so the subset
$\hh^s\times\cc^n$ is invariant under the action
of $\PPPP$.
This action is properly discontinuous
and the quotient is a complex analytic
manifold of dimension $s+n$
which will be denoted by 
$X(K,\OOOO, U)$
.
The notation $X(K,U)$ used in the article \cite{OT}
pertains to the case when the order $\OOOO$ 
is the maximal order of $K$.

\subsection{OT-manifolds  as manifolds of type $T(M,\DD)$}
\lb{su:ot-as-tm}

In this section we 
assume that there is a Dirichlet unit $\xi$ in $K$ 
\sut~ $\qq(\xi)=K$.
This assumption 
holds for example 
when $s+n\geq 2$ and there are no proper subfields 
$\qq\subsetneq K'\subsetneq K$; this is always the 
case if $s+n\geq 2$ and $2n+s$ is a prime number.
Replacing $\xi$ by $\xi^2$ if necessary
we can assume that $\xi$ is positive.
Denote by $P$ the minimal polynomial of $\xi$,
let $C_P$ be the companion matrix of $P$,
and put $B_P=C_P^\top$.
Denote by $\OOOO$ the order $\zz[\xi]$, and let $U$ be any subgroup 
of the group of units of $\zz[\xi]$,
isomorphic to $\zz^s$, and  generated by positive units $u_1, \ldots, u_s$.

 Choose polynomials $D_i\in \zz[t]$,
\sut~ $u_i=D_i(\xi)$. 
The family $\DD=(D_1, \ldots, D_s)$ is a Dirichlet family
for the matrix $C_P$
(see Proposition \ref{p:diri}
and Corollary \ref{c:irreducible-case}).

Both $X(K,\OOOO, U)$
and 
$T(M,\DD)$
are obtained as quotients
of actions of semidirect products of groups $\zz^s$ and $\zz^{s+2n}$
on $\hh^s\times \cc^n$.
We are now going to describe these two actions
in details and show that they are isomorphic.
Let $\a_1, \ldots, \a_s, \b_1, \ldots, \b_n,\bar\b_1, \ldots, \bar\b_n$
be the roots of $P$.
Renumbering the roots of $P$ if necessary
we can assume that 
$\s_i(\xi)=\a_i$ for $i\leq s$,
and
$\s_i(\xi)=\b_i, \ 
\s_{i+s}(\xi)=\ove{\b_i}
$ for $i> s$,

\pa
{\bf 1. }  \  The \ma~ $T(B_P, \DD).$
\pa
\noindent

The eigenvalues of the matrix $B_P$
are the same as those  of the matrix $C_P$, that is,
$\a_1, \ldots, \a_s,  \b_1, \ldots, \b_n, \bar\b_1, \ldots, \bar\b_n,$.
The corresponding eigenvectors of $B_P$ are:
$$
a_r=(1,\a_r,\ldots, \a_r^{2n+s-1}), \ \ 
b_i=(1,\b_i,\ldots, \b_i^{2n+s-1})
$$
(where $1 \leq r \leq s$ and 
$1 \leq i \leq n$).
We have therefore a decomposition
of $\cc^s\times\cc^{2n}$ into 
a direct sum of invariant 
subspaces
$$
\cc^{s+2n} = V\oplus W\oplus \bar W
$$
where $V$ is generated by the eigenvectors 
$
a_1, \ldots , a_s
$ 
 and $W$ is 
 generated by the eigenvectors 
$
b_1, \ldots , b_n.
$ 
The matrix $R$ is therefore diagonal 
with diagonal entries $\b_1, \ldots, \b_n$.
The vectors $v_i$
(see the formula \rrf{f:def-v})
are 
given by the formulas
$$
v_1=(1, \ldots, 1),
\
v_2=(\a_1, \ldots, \a_s, \b_1, \ldots, \b_n), 
\ 
\ldots 
\ 
,
$$
$$
v_{2n+s}=(\a_1^{2n+s-1},  \ldots, \a_s^{2n+s-1}, \  \b_1^{2n+s-1}, \ldots, \b_n^{2n+s-1}).
$$
\noindent
%Both  matrices $B_{P}, \ C_P$
%restricted to the subspace $W$ generated by 
%the vectors $v_i$ with $s+1\leq i \leq n+1$
%are diagonal in this basis
%and the diagonal entries are equal to 
%$\b_1, \ldots, \b_n$.
The element $g_0$ 
(see the formula \rrf{f:def-g0})
acts as follows:
$$
g_0\cdot (w_1, \ldots ,w_s,z_1, \ldots, z_n)
=
(\a_1w_1, \ldots \a_s w_s , \b_1z_1,\ldots,  \b_n z_n).
$$
The elements $g_i$ (where $1\leq i\leq s$)
act as follows:
%The matrix $D_i(B_P)$ acts as $D_i(g_0)$ (here $1\leq i\leq s$),
%that is, the action is by the diagonal 
%matrix with diagonal entries
$$
g_i\cdot (w_1, \ldots ,w_s,z_1, \ldots, z_n)
=
$$
$$
\Big(D_i(\a_1)w_1, \ldots ,D_i(\a_s) w_s , D_i(\b_1)z_1,\ldots, 
D_i(\b_n) z_n\Big).
$$

\pa
{\bf 2. }  \  The \ma~ $X(K,\OOOO, U).$
\pa
\noindent
The $\s$-image of $\xi^r$
equals
$(\a_1^r,\ldots, \a^r_s, \b^r_1, \ldots, \b^r_n)$;
 this vector  equals  $v_{r+1}$
(here $0\leq r\leq s+2n-1$).
Thus the group
of
translations
of
$\hh^s\times\cc^n$ 
by elements of 
$
\s(\zz[\xi])
$
is
generated
by the vectors $v_i$.
 The action of $\xi\in \OOOO$ on 
$\hh^s\times \cc^n$ is given by
$$
\xi\cdot (w_1, \ldots, w_s, z_1, \ldots, z_n)
$$
$$
=
\big(\s_1(\xi)w_1, \ldots, \s_s(\xi)w_s,\ \s_{s+1}(\xi)z_1, \ldots, \s_{n+s}(\xi)z_n\big)
$$
$$
=
(\a_1w_1, \ldots, \a_sw_s, \b_{1}z_1, \ldots, \b_{n}z_n).
$$
Therefore the elements $D_i(\xi)$ generating the group $U$ 
act as follows:
$$
D_i(\xi)(w,z)
=
\Big(D_i(\a_1)w_1, \ldots, D_i(\a_s)w_s, D_i(\b_1)z_1, \ldots, D_i(\b_n)z_s\Big).
$$
We deduce that the semidirect products 
of groups $\zz^s$ and $\zz^{2n+s}$
arising in the cases 1 and 2 are isomorphic 
and their actions on 
$\hh^s\times\cc^n$ are identical. 
Thus we arrive at the following 
conclusion.

\bepr\lb{p:ot-tm}
We have a biholomorphism
$$T(M,B_P)\approx X(K,\OOOO, U).\ \ \ \ \ \qs $$
\enpr

\subsection{The  non-diagonalizable case}
\lb{su:non-diag-bis}

\bele\lb{l:OT-diag-type}
Let $K$ be an algebraic number field
admitting $s$ real
and $n$
complex
embeddings into $\cc$  so that 
 $(K:\qq)=2n+s$.
Let $\OOOO$ be an order in $K$, denote by $\OOOO^*$
the multiplicative group of units of $\OOOO$.
Let  $U\approx \zz^s$ 
be a subgroup \sut~ every $\g\in U$ is positive.
Let $X=X(K,\OOOO, U)$ be the corresponding OT-manifold.
Then the group $\pi_1(X)$ is a semidirect product 
$\zz^s\ltimes\zz^{2n+s}$ of diagonal type.
\enle
\Prf
The group $\pi_1(X)$ is 
isomorphic to $G_{s,2n}(A)$ where $A:U\to GL(2n+s,\zz)$
denotes the action of the group $U$ on the 
order $\OOOO\approx \zz^{2n+s}$.
Let $\g\in U$, and $P$ be the minimal polynomial of $\g$.
The roots of $P$ are simple and different from 1. 
Since $P(A(\g))=0$, the minimal polynomial of 
the matrix $A(\g)$ has the same properties.
Therefore $A(\g)$ is diagonalizable  and Lemma is proved. $\qs$

The following proposition is now a direct consequence of 
Proposition \ref{p:isos}.

\bepr\lb{p:compar}
Let $M\in SL(2n+s,\zz)$.
Assume that there is a special Dirichlet
family 
$\DD=(D_1, \ldots , D_s)$
for $M$
\sut~ at least one of matrices
$D_i(M)$ 
is 
non-diagonalizable over $\cc$.
Then the group $\pi_1(T(M,\DD))$
is not isomorphic to 
the fundamental group of
any of manifolds
$X(K,\OOOO,U)$ constructed in \cite{OT}.
Therefore $T(M,\DD)$ is not homeomorphic 
to any of manifolds
$X(K,\OOOO,U)$.$\qs$
\enpr
The assumptions of the previous proposition 
hold if $\DD$ is primary. Therefore we obtain the next Corollary.

\beco\lb{c:compar2}
Let $M\in SL(2n+s,\zz)$ be a 
non-diagonalizable matrix of type $\JJ$,
and $\DD$  a primary Dirichlet family for $M$.
Then $T(M,\DD)$ is not homeomorphic 
to any of manifolds
$X(K,\OOOO,U)$.$\qs$
\enco

\section{Acknowledgements}
\lb{s:ack}

The first author thanks the Nantes University and the 
DefiMaths program for the support and warm hospitality. 
The first author was partially supported by JSPS KAKENHI
Grant Numbers 16K05142, 17H06128, 26287013, 19H01788. 
The second author thanks Professor F. Bogomolov for 
initiating him to the theory of Inoue surfaces in 2013,
for several discussions on this subject and for support.
The work on this article began in January 2018 
when the second author was visiting the Tokyo Institute of
Technology; many thanks for the warm hospitality and support.
We are grateful to K. Conrad and to Th. Laffey, whose 
works \cite{conrad1}, \cite{Laffey}
helped a lot  with the algebraic 
part of our project.

\end{document}